\documentstyle{article}
\begin{document}
\baselineskip+6pt \small \begin{center}{\textit{ In the name of
Allah, the Beneficent, the Merciful}}\end{center} \large

\begin{center}{ON THE FIELD OF DIFFERENTIAL RATIONAL INVARIANTS OF A SUBGROUP
OF AFFINE GROUP(PARTIAL DIFFERENTIAL CASE)} \end{center}

\begin{center} {Ural Bekbaev \footnote[1]{e-mail: bekbaev@science.upm.edu.my}}
\end{center}

\begin{center} {Department of Mathematics $\&$ Institute for Mathematical Research,}\end{center}
\begin{center} {FS, UPM, 43400, Serdang, Selangor, Malaysia.}\end{center}

\begin{abstract}{An differential field  $(F;\partial_1,...,\partial_m)$ of characteristic zero, a
subgroup $H$ of affine group $ GL(n,C)\propto C^n$ with respect to
its identical representation in $F^n$ and the following two fields
of differential rational functions in $x=(x_1,x_2,...,x_n)$-column
vector,
$$C\langle x, \partial \rangle^H=\{f^{\partial}\langle x\rangle \in C\langle x, \partial\rangle :
f^{\partial}\langle hx+ h_0\rangle = f^{\partial}\langle x\rangle\
\mbox{whenever} \ (h,h_0)\in H \},$$
$$C\langle x, \partial\rangle^{(GL^{\partial}(m,F),H)}=\{f^{\partial}\langle x\rangle \in C\langle x, \partial\rangle :
f^{g^{-1}\partial}\langle hx+ h_0\rangle = f^{\partial}\langle
x\rangle\ \mbox{whenever} \ g\in GL^{\partial}(m,F) \ \mbox{and} \
(h,h_0)\in H \}$$ are considered, where $C$ is the constant field
of $(F,\partial)$, $C\langle x, \partial\rangle$ is the field of
$\partial$-differential rational functions in $x_1,x_2,...,x_n$
over $C$ and $$GL^{\partial}(m,F)= \{
g=(g_{jk})_{j,k=\overline{1,m}}\in GL(m,F):
\partial_ig_{jk}= \partial_jg_{ik}\ \mbox{for} \
i,j,k=\overline{1,m}\}$$, $\partial$ stands for the column-vector
with the "coordinates" $\partial_1,. . .,\partial_m$. The field
$C\langle x, \partial\rangle^H$ ($C\langle x,
\partial\rangle^{(GL^{\partial}(m,F),H)}$) is an important tool in the
equivalence problem of patches( respect. surfaces) in Differential
Geometry with respect to the motion group $H$. In this paper a
pure algebraic approach is offered to describe these fields. The
field $C\langle x, \partial\rangle^{(GL^{\partial}(m,F),H)}$ is
considered and investigated as a differential field with respect
to a commuting system of differential operators
$\delta_1,...,\delta_m$. Its relation with differential field
$(C\langle x,
\partial\rangle^H,\partial)$ is shown. It is shown also that
$C\langle x,
\partial\rangle^H$ can be derived from some algebraic ( without
derivatives) invariants of $H$. \vspace{0.1cm}
} \end{abstract}
\vspace{0.5cm}

{\bf  1. Introduction.}

Let $n$, $m$  be natural numbers and $H$ be a subgroup of affine
group \\ $ GL(n,R)\propto R^{n}$, $G= Diff(B)$ be the group of
diffeomorphisms of the open unit ball $B\subset R^{m}$,
$u:B\rightarrow R^{n}$ be a surface, where $u$ is considered to be
infinitely smooth.

A function $f^{\partial}(u(t))$ of $u(t)=(u_1(t),. . ., u_{n}(t))$
and its finite number of derivatives relative to $\partial_1=
\frac{\partial}{\partial t_1}, . .
.,\partial_m=\frac{\partial}{\partial t_m}$ is said to be
invariant(more exactly, $ (G,H)$- invariant) if the equality
$$f^{\partial}(u(t))= f^{\delta}(hu(s(t)) + h_o)$$ is valid for
any $s\in G$, $(h,h_o)\in H$ and $t\in B$,
where $u(t)$ stands for the column vector with coordinates
$u_1(t),. . ., u_{n}(t))$, $ s(t)= (s_1(t),. . .,s_m(t))$,
$\delta_i= \frac{\partial}{\partial s_i}$.

Let $t$ run $B$ and $F= C^{\infty}(B, R)$ be the differential ring
of infinitely smooth functions relative to differential operators
$\partial_1= \frac{\partial}{\partial t_1}, . .
.,\partial_m=\frac{\partial}{\partial t_m}$. The constant ring of
this differential ring is $R$ i.e.  $$R= \{a\in F: \partial_i a=
0\quad \mbox{at}\quad i=\overline{1,m} \}$$ Every infinitely
smooth surface $u:B\rightarrow R^{n}$ can be considered as an
element of differential module $(F^{n};\partial_1, \partial_2,. .
.,\partial_m)$, where $\partial_i= \frac{\partial}{\partial t_i}$
acts on elements of $F^{n}$ coordinate-wisely. If elements of this
module are considered as column vectors the above transformations,
involved in definition of invariant function, look like $u=(u_1,.
. .,u_{n})\mapsto hu+ h_0,\quad \partial \mapsto g^{-1}\partial$
as far as $$\frac{\partial}{\partial t_i}=
\sum_{j=1}^{m}\frac{\partial s_j(t)}{\partial
t_i}\frac{\partial}{\partial s_j(t)}$$ , where $g$ is matrix with
the elements  $g_{ij}=\frac{\partial s_j(t)}{\partial t_i}$
 at $i,j=\overline{1,m}$, $\partial$  is the column vector with the "coordinates" $\frac{\partial}{\partial
t_1},. . .,\frac{\partial}{\partial t_m}$ , $(h,h_0)\in H$, $s\in
G$. Moreover $\partial_ig_{jk}= \partial_jg_{ik}\ \mbox{for} \
i,j,k =\overline{1,m}.$

Therefore for any differential field $(F;\partial_1, \partial_2,.
. .,\partial_m)$, i.e. $F$ is a field and $\partial_1,
\partial_2,. . .,\partial_m$ is a given commuting with each other
system of differential operators on $F$, one can consider the
transformations
$$u=(u_1,. . .,u_{n})\mapsto hu+ h_0,\qquad \partial
\mapsto g^{-1}\partial,$$ , where $u \in F^n,$ $(h,h_0)\in H$ is a
given subgroup of affine group $ GL(n,C)\propto C^{n}$, $g\in
GL^{\partial}(m,F)$, \[GL^{\partial}(m,F)= \{ g\in  GL(m,F):
\partial_ig_{jk}= \partial_jg_{ik}\ \mbox{for} \ i,j,k=\overline{1,m}\}\]
$\partial$ stands for the column-vector with the "coordinates"
$\partial_1,. . .,\partial_m$, $C$ is the constant field of
$(F;\partial )$ i.e.
$$C= \{a\in F: \partial_i a= 0\quad \mbox{at}\quad
i=\overline{1,m} \}.$$

It should be noted that for any $g \in GL^{\partial}(m,F)$ the
differential operators $\delta_1,\delta_2,. . .,\delta_m$, where
$\delta= g^{-1}\partial$, also commute with each other. So for any
$g\in GL^{\partial}(m,F)$ one can consider the differential field
$(F,\delta)$, where $\delta= g^{-1}\partial$. This transformation
is an  analogue of gauge transformations for abstract differential
field $(F,\partial)$.

{\bf Remark 1.} In common case the set $GL^{\partial}(m,F)$ is not
a group with respect to the ordinary product of matrices as far as
it is not closed with respect to that product. But by the use of
it a natural groupoid ([1]) can be constructed with the base
$\{g^{-1}\partial :  g\in GL^{\partial}(m,F)\}$.

{\bf Remark 2.} Let $g \in GL^{\partial}(m,F)$ and $\delta=
g^{-1}\partial$. It is clear that $$\{a\in F: \partial_1
a=...=\partial_m a=0 \}=C=\{a\in F: \delta_1 a=...=\delta_m
a=0\}$$ One interesting question is: When does one have the
equality
$$ \bigcup_{k\in N}\{a\in F: \partial^{\alpha}a=0 \ \mbox{whenever}\ \vert \alpha \vert=k
\}=\bigcup_{k\in N}\{a\in F: \delta^{\alpha}a=0 \ \mbox{whenever}\
\vert \alpha \vert=k \} \ ?$$, where $\alpha
=(\alpha_1,\alpha_2,...,\alpha_m)$, $\vert \alpha \vert
=\alpha_1+\alpha_2+...\alpha_m$, $\alpha_i$ are nonnegative
integers and
$\delta^{\alpha}=\delta_1^{\alpha_1}\delta_2^{\alpha_2}...\delta_m^{\alpha_m}.$

 For example, if
$g \in GL^{\partial}(m,F)$, $\det g\in C$ and all entries of $g$
are in $$ \bigcup_{k\in N}\{a\in F: \partial^{\alpha}a=0 \
\mbox{whenever}\ \vert \alpha \vert=k \}$$ then does it imply the
above equality? Of course it is an generalization of the famous
Jacobian Conjecture which is the same problem when
$F=Q(t_1,...,t_m)$, $\partial_1=\frac{\partial}{\partial
t_1}$,..., $\partial_m=\frac{\partial}{\partial t_m}$.

Let in future $x_1,. . .,x_{n}$ be differential algebraic
independent variables over $F$ and $x$ stand for the column vector
with coordinates $x_1,. . .,x_{n}$. We use the following notations
: $C[x]$ - the ring of polynomials in $x_1,. . .,x_{n}$ (over
$C$), $C(x)$- the field of rational functions in $x$, $C\{
x,\partial \}$ -the ring of $\partial$-differential polynomial
functions in $x$ and $C\langle x,\partial \rangle$-is the field of
$\partial$-differential rational functions in $x$ over $C$.

{\bf Definition.} \textit{ An element $f^{\partial}\langle
x\rangle \in C\langle x, \partial\rangle$ is said to be
 $(GL^{\partial}(m,F),H)$- invariant \\ $(GL^{\partial}(m,F)$- invariant; $H$- invariant)
 if the equality
$$f^{g{-1}\partial}\langle hx+h_0\rangle = f^{\partial}\langle
x\rangle $$ (respect. $f^{g{-1}\partial}\langle x\rangle =
f^{\partial}\langle x\rangle; f^{\partial}\langle hx+h_0\rangle =
f^{\partial}\langle x\rangle )$ is valid for any $g\in
GL^{\partial}(m,F),\quad (h,h_0)\in H$.}

Let $C\langle x,
\partial\rangle^{(GL^{\partial}(m,F),H)}$  ($C\langle x, \partial\rangle^{GL^{\partial}(m,F)}$,
$C\langle x, \partial \rangle^H$ ) stand for the set of all such
$(GL^{\partial}(m,F),H)$- invariant (respect.
$GL^{\partial}(m,F)$- invariant, $H$- invariant) elements of $
C\langle x,\partial\rangle$.

The fields $C\langle x,
\partial\rangle^{(GL^{\partial}(m,F),H)}$, $C\langle x, \partial \rangle^H$
and their relations
are investigated in [2] for the case of $m=1$. Some results on
these fields can be found in [3] for the case of $m=n-1$. In the
case of finite $H$ more strong results than results of this paper
are presented in [4]. The first variant of this paper is published
in [5]. The needed notions and results from Differential Algebra
can be found in [6].

{\bf  2. Preliminary}

 In future let
$(F,\partial )$ stand for a field $F$ with fixed commuting system
of differential operators $\partial_1,. . .,\partial_m$ and $C$ be
its constant field i.e. $C=\{ a\in F:\partial_1a=. .
.=\partial_ma=0\}$.

{\bf Proposition 1.} \textit{If the system of differential
operators $\partial_1,. . .,\partial_m$ is linear independent over
$F$ then the differential operators $\delta_1,. .
.,\delta_m$,where $\delta= g^{-1}\partial, g\in GL(m,F)$, commute
with each other if and only if $  g\in GL^{\partial}(m,F)$.}

{\bf Proof.} Let $g\in GL(m,F)$,\ $\delta= g^{-1}\partial $. It is
clear that linear independence of $\partial_1,. . .,\partial_m$
implies linear independence of $\delta_1,. . .,\delta_m.$ For any
$i,j=\overline{1,m}$ we have $\partial_j=
\sum_{k=1}^mg_{jk}\delta_k$,\ $\partial_i\partial_j=
\sum_{k=1}^m(\partial_i(g_{jk})\delta_k+
g_{jk}\partial_i\delta_k)= \sum_{k=1}^m\partial_i(g_{jk})\delta_k
+ \sum_{k=1}^m\sum_{s=1}^mg_{jk}g_{is}\delta_s\delta_k$. Therefore
due to $\partial_i\partial_j= \partial_j\partial_i $one has
\begin{eqnarray}\begin{array}{c}
\sum_{k=1}^m\partial_i(g_{jk})\delta_k +
\sum_{k=1}^m\sum_{s=1}^mg_{jk}g_{is}\delta_s\delta_k=
\sum_{k=1}^m\partial_j(g_{ik})\delta_k +
\sum_{k=1}^m\sum_{s=1}^mg_{jk}g_{is}\delta_k\delta_s
\end{array}\end{eqnarray}

If $\delta_k\delta_s = \delta_s\delta_k$ for any $k,s=
\overline{1,m}$ them due (1) one has
$\sum_{k=1}^m\partial_i(g_{jk})\delta_k  =
\sum_{k=1}^m\partial_j(g_{ik})\delta_k $ i.e. $\partial_i(g_{jk})=
\partial_j(g_{ik})$ for any  $i,j,k= \overline{1,m}$ because of linear independence of
$\delta_1,. . .,\delta_m$. Thus in this case
 $  g\in GL^{\partial}(m,F)$.

Vice versa, if $\partial_i(g_{jk})=
\partial_j(g_{ik})$ for any  $i,j,k= \overline{1,m}$  then due (1) at any $a\in F$ one has
$\sum_{k=1}^m\sum_{s=1}^mg_{jk}g_{is}\delta_s\delta_k a=
\sum_{k=1}^m\sum_{s=1}^mg_{jk}g_{is}\delta_k\delta_s a$ for any
$i,j= \overline{1,m}$. These equalities can be written in the
following matrix form $g(\delta_1\delta a,. . .,\delta_m\delta
a)g^t= g(\delta_1\delta a,. . .,\delta_m\delta a)^tg^t $, where
$t$ means transposition. Therefore $(\delta_1\delta a,. .
.,\delta_m\delta a)= (\delta_1\delta a,. . .,\delta_m\delta a)^t$
i.e.
 $\delta_k\delta_s a= \delta_s\delta_k a$ for any
$k,s= \overline{1,m}$, which completes the proof of Proposition 1.

It should be noted that for $g\in GL^{\partial}(m,F)$ and $\delta=
g^{-1}\partial$ the following equality is valid
\begin{eqnarray}\begin{array}{c} GL^{\delta}(m,F)=
g^{-1}GL^{\partial}(m,F). \end{array}\end{eqnarray}

If $(K,d)$ is an ordinary differential field of characteristic
zero with a constant field $$K_0= \{a\in K: d(a)=0\}$$ then the
following criterion  is well known: A system of elements
$b_1,b_2,...,b_n $ of $K$ is $K_0$-linear dependent if and only if
$$\det [b,d(b),...,d^{n-1}(b)]=0$$, where $b$ stands for the
vector $(b_1,b_2,...,b_n) $. Similar question can be asked in
common case: If $b_1,b_2,...,b_n $ is a system of elements of a
differential field $(F,\partial )$ how one can find out if it is
linear dependent over $C$? The following result deals with this
problem in common case.

Consider  the differential field $(F,\partial )$ of characteristic
zero, its constant field $C$ and indeterminates $\{d^it_j: i\in W,
j=\overline{1,m}\}$. Let $F\langle t \rangle $ ($F\{t\}$) stand
for the field(respect. ring) of all rational(respect.
polynomial)functions in  $\{d^it_j: i\in W, j=\overline{1,m}\}$
over $F$. One can make it an ordinary differential field(respect.
ring) $F\langle t,d \rangle $ (respect. $F\{t,d\}$) by allowing

1. $d(a)= \sum_{i=1}^m\partial_i(a)dt_i $ for any $a\in F$.

2. $d(d^it_j)= d^{i+1}t_j$ for any $ i\in W, j=\overline{1,m}$.

The following result shows that the introduction of the ordinary
differential field $F\langle t,d \rangle $ reduces the above
stated question once again to the ordinary case.

{\bf Proposition 2}. \textit{ The constant field of the ordinary
differential field $F\langle t,d \rangle $ is the came $C$.}

 This result can be deduced easily from the fact that
$f^d\{t\}$ divides $df^d\{t\}$ if and only if $f^d\{t\}\in F $,
where $f^d\{t\}\in F\{t\}$. As far as $F\subset F\{t\}$ the answer
to the above question can be given in the following way: The
system  $b_1,b_2,...,b_n $ of elements $F$ is $C$-linear dependent
if and only if $\det [b,d(b),...,d^{n-1}(b)]=0$.

{\bf Proposition 3}. \textit{Let  $ (F,\partial_1,
\partial_2, ..., \partial_m)$- be a differential field of characteristic zero. The following three properties are equivalent.\\
a) The system of differential operators  $ \partial_1, \partial_2,
 ...,\partial_m $ is linear independent over $ F$.\\
b) There is no nonzero $ \partial $- differential
 polynomial over  $ F$ which vanishes at all values
 of indeterminates from  $ F$.\\
c) If $p^{\partial}\{x_{11}, x_{12},..., x_{1m},
 x_{21},..., x_{2m},..., x_{m1},..., x_{mm}\}=
p^{\partial}\{(x_{ij})_{i,j= \overline{1,m}}\}$ is a differential
polynomial over $ F$  such that $ p^{\partial}\{g\}= 0 $  at any
$ g= (g_{ij})_{i,j= \overline{1,m}}\in GL^{\partial}(m,F)$, then
  $ p^{\partial}\{(t_{ij})_{i,j=
 \overline{1,m}}\}= 0 $ is also valid,
for any indeterminates  $ (t_{ij})_{i,j= \overline{1,m}}$,  for
which $ \partial_k t_{ij}=
 \partial_i t_{kj}$  at  $ i,j,k= \overline{1,m}.$ }

{\bf Proof.} The equivalence of properties a) and b) is proved in
[6, p.139].

It is evident that c) implies b). In fact if there are nonzero $
\partial $- differential
 polynomials over  $ F$ which vanish at all values
 of indeterminates from  $ F$ we can take one with
minimal number of variables. Let $f\{z_1, z_2,. . .,z_l\}$ be such
a polynomial. If $l> 1$ and $f\{z_1, a_2,. . .,a_l\} \neq 0$ for
some $a_2,. . .,a_l\in F$ then it contradicts minimality of $l$.
Because the nonzero polynomial in one variable $f\{z_1, a_2,. .
.,a_l\}$ will vanish at all values of $z_1$ from $F$. If $l> 1$
and $f\{z_1, a_2,. . .,a_l\}= 0$ for all $a_2,. . .,a_l\in F$ then
considering $f\{z_1, z_2,. . .,z_l\}$ as a $ \partial $-
differential polynomial in $z_1$ over $F\{z_2, z_3,. . .,z_l\}$
once again we will have a contradiction. Indeed in this case at
least one of the coefficients of this polynomial has to be nonzero
$ \partial $- differential polynomial in $z_2, z_3,. . .,z_l$ ( as
$f\{z_1, z_2,. . .,z_l\}$ is a nonzero polynomial)and vanish at
all values of $z_2, z_3,. . .,z_l$ from $F$. It contradicts
minimality of $l$. Thus $l= 1$ and we can consider nonzero
polynomial $p^{\partial}\{(x_{ij})_{i,j= \overline{1,m}}\}=
f\{x_{11}\}$ which vanishes at any $ g= (g_{ij})_{i,j=
\overline{1,m}}$  $\in GL^{\partial}(m,F)$. This contradicts
property c).

Let us prove now that b) implies c). Assume that $
p^{\partial}\{(t_{ij})_{i,j= \overline{1,m}}\}\neq 0 $ for some
polynomial $p^{\partial}\{(x_{ij})_{i,j= \overline{1,m}}\}$. Due
to the equalities $ \partial_{k}t_{ij}=
 \partial_{i}t_{kj}$, $i,j,k= \overline{1,m} $
the nonzero $ p^{\partial}\{(t_{ij})_{i,j= \overline{1,m}}\} $ can
be represented as a polynomial $P$ of the monomials $
 {\partial_k}^{n_{k,i}} {\partial_{k+1}}^{n_{k+1,i}} ...
{\partial_m}^{n_{m,i}} t_{ki}, $ where  $ n_{j,i}$- are
nonnegative integers, $i,k= \overline{1,m} $.  Let $t_1, t_2,...,
t_m $ be any differential indeterminates over $ F$.  The
 inequality $ 0\neq p^{\partial}\{(t_{ij})_{i,j=
 \overline{1,m}}\} $ and substitution $ t_{ij}=
 \partial_it_j $ give us a nonzero differential polynomial $\det(\partial_it_j)_{i,j=
 \overline{1,m}}P$
 in $t_1, t_2,..., t_m  $ the value of which at any
 $(a_1, a_2,..., a_m)  $ from $F^m$ is zero. This contradicts
 property b).

{\bf Proposition 4}. \textit{If  $a=(a_1, a_2,..., a_m)$ and
$b=(b_1, b_2,..., b_m)$ are any two nonzero row vectors from $F^m$
then there is such extension $(F_1,\partial)$ of $(F,\partial)$
where the equation $aT=b$ has solution in $GL^{\partial}(m,F_1)$.}

{\bf Proof.} Assume , for example, that $a_1\neq 0$ and $
\{t_{ij}\}_{i= \overline{2,m}, j= \overline{1,m}}$ are such
 differential indeterminates over $ F$ that $
 \partial_{k}t_{ij}= \partial_{i}t_{kj}$  for  $i,k=
 \overline{2,m}, j= \overline{1,m}$.

Consider $$(a_1, a_2,..., a_m)\pmatrix{y_1& y_2&.&.&.&y_m\cr
t_{21}& t_{22}&.&.&.&t_{im}\cr .& .&.&.&.&.\cr t_{m1}&
t_{m2}&.&.&.&t_{mm}\cr}=(b_1, b_2,..., b_m)$$ as a system of
linear equations in $y_1, y_2,..., y_m$.
 It has solution
\[(y_1, y_2,..., y_m)=(t_{11}, t_{12},. . ., t_{1m})= \frac{1}{a_1}(b- \sum_{i=2}^{m}a_i(t_{i1}, t_{i2},.
. ., t_{im})) \]and the determinant of the corresponding matrix
$T=(t_{ij})_{i,j=\overline{1,m}}$ is equal to
$$\frac{1}{a_1}\det\pmatrix{b_1& b_2&.&.&.&b_m\cr t_{21}&
t_{22}&.&.&.&t_{im}\cr .& .&.&.&.&.\cr t_{m1}&
t_{m2}&.&.&.&t_{mm}\cr}$$ which is not zero because of  $b\neq 0$.
Furthermore if one defines $\partial_1(t_{k1}, t_{k2},. . .,
t_{km})$ as
\[\partial_1(t_{k1}, t_{k2},. . ., t_{km})= \partial_k (t_{11}, t_{12},. . ., t_{1m})= \partial_k(\frac{1}{a_1}(b-\sum_{i=2}^{m}a_i(t_{i1}, t_{i2},.
. ., t_{im}))) \] then $T\in GL^{\partial}(m, F_1)$, where
$F_1=F\langle
\{t_{i,j}\}_{i=\overline{2,m},j=\overline{1,m}};\partial \rangle$.
This is the proof of Proposition 4.

In future let $e$ stand for the row vector $(1,0,0,...,0)\in F^m$
and $T$ stand for the matrix $(t_{ij})_{i,j= \overline{1,m}}$,
where $\{t_{ij}\}_{i,j= \overline{1,m}}$- are differential
indeterminates with the basic relations $ \partial_{k}t_{ij}=
 \partial_{i}t_{kj}$  for all $i,j,k= \overline{1,m}.$

{\bf Corollary}. \textit{If $ (F; {\partial}_1,
 {\partial}_2, ..., {\partial}_m) $- is a
 differential field of characteristic zero,
  $ \partial_1, \partial_2, ..., \partial_m $
 is linear independent over  $ F$ and $
 p^{\partial}\{t_1, t_2,..., t_m, t_{m+1},...,t_{m+l}\} $- is an arbitrary
 nonzero differential polynomial over $ F$ then} $\quad a) \  p^{\partial}\{eT,t_{m+1},..., t_{m+l} \} \neq 0,$

 \[b)\ p^{T^{-1}\partial}\{eT^{-1}, t_{m+1},..., t_{m+l} \} \neq 0 ,\qquad  c)\ p^{T^{-1}\partial}\{eT,t_{m+1},..., t_{m+l}\} \neq 0.  \]

{\bf Proof.} The proof of inequality a) is evident due to
Propositions 3 and 4.

 Let us prove b). If one assumes that
$ p^{T^{-1}\partial}\{eT^{-1}, t_{m+1},..., t_{m+l}\}=0$ then in
particular for $T_0=(\partial_i x_j)_{i,j=\overline{1,m}}$ one has
$ p^{T^{-1}_0\partial}\{eT^{-1}_0, t_{m+1},..., t_{m+l}\}=0$. Once
again it will have to remain be true if one substitutes
$T^{-1}\partial$ for $\partial$ into it. But $\delta_0=
T_0^{-1}\partial$ is invariant with respect to such substitution
and $T_0$ is transformed to $T^{-1}T_0$ so $$
p^{T^{-1}_0\partial}\{eT^{-1}_0T, t_{m+1},..., t_{m+l}\}=0.$$ But
for any $S_0\in GL^{\delta_0}(m, F\langle x_1,...,x_m;\partial
\rangle )$ the equation $T^{-1}_0T= S_0$ has solution in \\
$GL^{\partial}(m, F\langle x_1,...,x_m;\partial \rangle )$, namely
$T= T_0S_0$. Therefore due to Proposition 3 for the matrix of
variables $S=(s_{ij})_{i,j=\overline{1,m}}$ for which
$\delta_{0i} s_{jk}= \delta_{0j} s_{ik}$ for all
$i,j,k=\overline{1,n}$ one has
$$p^{\delta_0}\{eS, t_{m+1},..., t_{m+l}\}=0.$$ Due to the Corollary, part a), $p^{\delta_0}\{t_1,..., t_{m+l}\}=0$ i.e. $p^{\partial}\{t_1,..., t_{m+l}\}=0$ which is a contradiction. The proof of c) can be done in a similar way.

{\bf  3. On $(GL^{\partial}(m,F),H)$- invariants.}

In future it is assumed that $(F,\partial )$, where $\partial
=(\partial_1,. . .,\partial_m)$, is such a differential field
that:

1. Characteristic of $F$ is zero.

2. The system $\partial_1,. . .,\partial_m$ is linear independent
over $F$.

We use the following obvious fact repeatedly: If $t_1,...,t_l$ is
a $\partial$-algebraic independent system of variables   over $F$,
$g\in GL^{\partial}(m, F)$ and $p^{\partial}\{t_1,...,t_l\}$ is a
$\partial$-polynomial over $F$ then the following equalities are
equivalent.
\[p^{\partial}\{t_1,...,t_l\}=0, \ p^{g^{-1}\partial}\{t_1,...,t_l\}=0,\ p^{T^{-1}\partial}\{t_1,...,t_l\}=0. \]

In future let us assume that for the given subgroup $H$ of $
GL(n,C)\propto C^{n}$ we have such a nonsingular matrix
\[\Phi^{\partial}\langle x\rangle = (\phi^{\partial}_{ij}\langle
x\rangle)_{i,j=\overline{1,m}}\] , where
$\phi^{\partial}_{ij}\langle x\rangle\in C\langle x,\partial
\rangle$ and
 $\partial_{k}\phi^{\partial}_{ij}= \partial_{i}\phi^{\partial}_{kj}$  for  $i,j,k= \overline{1,m}$, that  \begin{eqnarray}\begin{array}{c}
\Phi^{g^{-1}{\partial}}\langle hx+h_0\rangle =
g^{-1}\Phi^{\partial}\langle x\rangle
\end{array}\end{eqnarray}
for any $g\in GL^{\partial}(m,F)$ and $(h,h_0)\in H$.

{\bf Remark 3.} For the given $H$ the existence of the none
singular matrix $\Phi^{\partial}\langle x\rangle$ with property
(3) is another problem. The paper does not touch this problem.
Existence problem of such matrix is considered in [3] in the case
of $n=m+1.$

It is evident that  $(C\langle x,\partial\rangle^H,\partial)$ is a
finitely generated $\partial$-differential field over $C$ as a
subfield of $(C\langle x,\partial \rangle,\partial )$
 and $C\langle x,\partial\rangle^{(GL^{\partial}(m,F),H)}$ is a differential field with respect to $\delta = \Phi^{\partial}\langle x\rangle^{-1}\partial$. One of the most important questions is the differential-algebraic transcendence degree of $C\langle x,\partial\rangle^{(GL^{\partial}(m,F),H)}$ as a such $\delta $-field over $C$.

{\bf Theorem 1.} \textit{$\delta$-tr.deg.$C\langle x,
\partial\rangle^{(GL^{\partial}(m,F),H)}/C= n-m$}

{\bf Proof.} First of all let us show that the system
 $\phi^{\partial}_{11}\langle x\rangle ,\phi^{\partial}_{12}\langle x\rangle,...,\phi^{\partial}_{1m}\langle x\rangle $ is  $\delta$-algebraic independent over $C\langle x, \partial\rangle^{GL^{\partial}(m,F)}$. Indeed if $p^{\delta}\{t_1,...,t_m\}$ is such a $\delta$-polynomial over $C\langle x, \partial\rangle^{GL^{\partial}(m,F)}$ for which $$p^{\delta}\{\phi^{\partial}_{11}\langle x\rangle ,\phi^{\partial}_{12}\langle x\rangle,...,\phi^{\partial}_{1m}\langle x\rangle\}=0$$ then it will have to remain be true if one substitutes $g^{-1}\partial$ for $\partial$ into it. Therefore, as far as all coefficients of $p^{\delta}\{t_1,...,t_m\}$, as well as $\delta$, are invariant with respect to such substitutions and $\Phi^{g^{-1}{\partial}}\langle x\rangle = g^{-1}\Phi^{\partial}\langle x\rangle $ one has
$p^{\delta}\{eg^{-1}\Phi^{\partial}\langle x\rangle\}=0$ i.e.
$p^{\delta}\{eT^{-1}\Phi^{\partial}\langle x\rangle\}=0$. But for
any $S_0\in GL^{\delta}(m, F\langle x, \partial\rangle)$ the
equation $\Phi^{\partial}\langle x\rangle^{-1}T=S_0$ has solution
in $GL^{\partial}(m, F\langle x, \partial\rangle)$, namely
$T=\Phi^{\partial}\langle x\rangle S_0$. It implies that for the
matrix of variables $S= (s_{ij})_{i,j=\overline{1,m}}$, for which
$\delta_i s_{jk}=\delta_j s_{ik}$ ,$i,j,k=\overline{1,m}$, one has
$p^{\delta}\{eS^{-1}\}=0$. Due to Corollary, part a) one has
$p^{\delta}\{t_1,...,t_m\}=0$.

Now let $f^{\partial}_1\langle x\rangle,..., f^{\partial}_l\langle
x\rangle $ be any system of elements of $C\langle x,
\partial\rangle^{GL^{\partial}(m,F)}$. We show that the system
$$\phi^{\partial}_{11}\langle x\rangle ,\phi^{\partial}_{12}\langle
x\rangle,...,\phi^{\partial}_{1m}\langle x\rangle,
f^{\partial}_1\langle x\rangle,..., f^{\partial}_l\langle x\rangle
$$ is $\delta$-algebraic independent over $C$ if and only if it is
$\partial$-algebraic independent over $C$.

Indeed if this system is $\delta$-algebraic independent over $C$
and $p^{\partial}\{t_1,...,t_{m+l}\}$ is any polynomial over $C$
for which $$p^{\partial}\{e\Phi^{\partial}\langle x\rangle ,
f^{\partial}_1\langle x\rangle,..., f^{\partial}_l\langle
x\rangle\}=0$$ then it will have to remain be true if one
substitutes $g^{-1}\partial$ for $\partial$ into it, where $g\in
GL^{\partial}(m,F)$. It implies that
$$p^{T^{-1}\partial}\{eT^{-1}\Phi^{\partial}\langle x\rangle ,
f^{\partial}_1\langle x\rangle,..., f^{\partial}_l\langle
x\rangle\}=0$$ because $f^{\partial}_1\langle x\rangle,...,
f^{\partial}_l\langle x\rangle$ are invariant with respect to such
transformations. But  \\ $T^{-1}\partial = (\Phi^{\partial}\langle
x\rangle^{-1}T)^{-1}\Phi^{\partial}\langle x\rangle^{-1}\partial =
(\Phi^{\partial}\langle x\rangle^{-1}T)^{-1}\delta $ and for any
$S_0\in GL^{\delta}(m,F\langle x;\partial\rangle)$ the equation
$\Phi^{\partial}\langle x\rangle^{-1}T=S_0$ has solution in
$GL^{\partial}(m, F\langle x, \partial\rangle)$, namely
$T=\Phi^{\partial}\langle x\rangle S_0$. Therefore for the matrix
of variables $S= (s_{ij})_{i,j=\overline{1,m}}$, for which
$\delta_i s_{jk}=\delta_j s_{ik}$ ,$i,j,k=\overline{1,m}$, one has
$$p^{S^{-1}\delta}\{eS, f^{\partial}_1\langle x\rangle,..., f^{\partial}_l\langle x\rangle\}=0.$$
Due to our assumption $f^{\partial}_1\langle x\rangle,..., f^{\partial}_l\langle x\rangle$ is $\delta$-algebraic independent over $C$ and therefore according to Corollary, part b),  $p^{\delta}\{t_1,...,t_{m+l}\}=0$ i.e. $p^{\partial}\{t_1,...,t_{m+l}\}=0.$ So the system\\
$\phi^{\partial}_{11}\langle x\rangle ,\phi^{\partial}_{12}\langle
x\rangle,...,\phi^{\partial}_{1m}\langle x\rangle,
f^{\partial}_1\langle x\rangle,..., f^{\partial}_l\langle x\rangle
$ has to be $\partial$-algebraic independent over $C$.

Vise versa, let $\phi^{\partial}_{11}\langle x\rangle
,\phi^{\partial}_{12}\langle
x\rangle,...,\phi^{\partial}_{1m}\langle x\rangle,
f^{\partial}_1\langle x\rangle,..., f^{\partial}_l\langle x\rangle
$ be $\partial$-algebraic independent over $C$. In this case first
of all the system $f^{\partial}_1\langle x\rangle,...,
f^{\partial}_l\langle x\rangle $ is $\delta$-algebraic independent
over $C$. Indeed \\$f^{\partial}_1\langle x\rangle,...,
f^{\partial}_l\langle x\rangle $ is $\partial$-algebraic
independent over $C\langle \{\phi^{\partial}_{ij}\langle x\rangle
\}_{i,j=\overline{1,m}};\partial \rangle$ as far as
$\partial_k\phi^{\partial}_{ij}\langle
x\rangle=\partial_i\phi^{\partial}_{kj}\langle x\rangle$ for
$i,j,k=\overline{1,m}$ and $\phi^{\partial}_{11}\langle x\rangle
,\phi^{\partial}_{12}\langle
x\rangle,...,\phi^{\partial}_{1m}\langle x\rangle,
f^{\partial}_1\langle x\rangle,..., f^{\partial}_l\langle x\rangle
$ is $\partial$-algebraic independent over $C$. But every nonzero
$p^{\delta}\{t_1,...,t_l\}$ over $C\langle
\{\phi^{\partial}_{ij}\langle x\rangle
\}_{i,j=\overline{1,m}};\partial \rangle$ can be considered as a
nonzero $\partial$- polynomial over $C\langle
\{\phi^{\partial}_{ij}\langle x\rangle
\}_{i,j=\overline{1,m}};\partial \rangle$. Therefore supposition
$p^{\delta}\{f^{\partial}_1\langle x\rangle,...,
f^{\partial}_l\langle x\rangle\}=0$ leads to a contradiction that
$f^{\partial}_1\langle x\rangle,..., f^{\partial}_l\langle
x\rangle $ is $\partial$-algebraic independent over $C\langle
\{\phi^{\partial}_{ij}\langle x\rangle
\}_{i,j=\overline{1,m}};\partial \rangle$.

Let us assume that  for some polynomial $p^{\delta}\{t_1,...,t_{m+l}\}$ over $C$ one has \\
$p^{\delta}\{e\Phi^{\partial}\langle x\rangle
,f^{\partial}_1\langle x\rangle,..., f^{\partial}_l\langle
x\rangle \}=0$. It should remain be true if one substitutes
$g^{-1}\partial$ for $\partial$ into it, where $g\in
GL^{\partial}(m,F)$, which leads to
$p^{\delta}\{eT^{-1}\Phi^{\partial}\langle x\rangle
,f^{\partial}_1\langle x\rangle,..., f^{\partial}_l\langle
x\rangle \}=0$. But the equation $\Phi^{\partial}\langle
x\rangle^{-1}T=S_0$ has solution in $GL^{\partial}(m, F\langle x,
\partial\rangle)$ for any $S_0\in GL^{\delta}(m,F\langle
x;\partial\rangle)$ and therefore
$$p^{\delta}\{eS^{-1}, f^{\partial}_1\langle x\rangle,..., f^{\partial}_l\langle x\rangle\}=0.$$
Now take into consideration that $f^{\partial}_1\langle
x\rangle,..., f^{\partial}_l\langle x\rangle $ is
$\delta$-algebraic independent  over $C$ and Corollary, part a) to
see that $p^{\delta}\{t_1,...,t_{m+l}\}=0$. So it is shown that
the system $$\phi^{\partial}_{11}\langle x\rangle
,\phi^{\partial}_{12}\langle
x\rangle,...,\phi^{\partial}_{1m}\langle x\rangle,
f^{\partial}_1\langle x\rangle,..., f^{\partial}_l\langle x\rangle
$$ is $\delta$- algebraic independent over $C$ if and only if it
is $\partial$-algebraic independent over $C$. In particular it
shows that the existence of $\Phi^{\partial}\langle x\rangle$ with
property (3) implies that $m\leq n$, because already we have got
that the system $\phi^{\partial}_{11}\langle x\rangle
,\phi^{\partial}_{12}\langle
x\rangle,...,\phi^{\partial}_{1m}\langle x\rangle$ is $\delta$-
algebraic independent over $C$.

It is evident that the system of components of the matrix
$\Phi^{\partial}\langle x\rangle^{-1}=(\psi^{\partial}_{ij}\langle
x\rangle )_{i,j=\overline{1,m}}$ generates $C\langle x;\partial
\rangle$ over $C\langle x, \partial\rangle^{GL^{\partial}(m,F)}$
as a $\delta$-differential field and $\delta_k
\psi^{\partial}_{ij}\langle x\rangle =\delta_i
\psi^{\partial}_{kj}\langle x\rangle$ for all
$i,j,k=\overline{1,m} $, which implies that
$\delta$-tr.deg.$C\langle x;\partial \rangle /C\langle
x;\partial\rangle^{GL^{\partial}(m,F)} \leq m$. But it already has
been established that $\phi^{\partial}_{11}\langle x\rangle
,\phi^{\partial}_{12}\langle
x\rangle,...,\phi^{\partial}_{1m}\langle x\rangle$ is
$\delta$-algebraic independent over $C\langle
x;\partial\rangle^{GL^{\partial}(m,F)}$ and therefore in reality
$\delta$-tr.deg.$C\langle x;\partial \rangle /C\langle
x;\partial\rangle^{GL^{\partial}(m,F)}= m$.

As a $\partial$-differential field $C\langle x;\partial\rangle$
over $C$ is generated by the elements of $C\langle
x;\partial\rangle^{GL^{\partial}(m,F)}$, as far as $x_1,...,x_n$
belong to it, and $\partial$-tr.deg.$C\langle x;\partial\rangle
/C=n.$ It means that one can find such a system
$f^{\partial}_1\langle x\rangle,..., f^{\partial}_{n-m}\langle
x\rangle $ elements of $C\langle
x;\partial\rangle^{GL^{\partial}(m,F)}$ for which the system
$$\phi^{\partial}_{11}\langle x\rangle ,\phi^{\partial}_{12}\langle x\rangle,...,\phi^{\partial}_{1m}\langle x\rangle , f^{\partial}_1\langle x\rangle,..., f^{\partial}_{n-m}\langle x\rangle $$
is $\partial$-algebraic independent over $C$. As it has been shown
that in this case it is $\delta$-algebraic independent over $C$ as
well. Therefore $\delta$-tr.deg.$C\langle x;\partial\rangle /C =n$
and $\delta$-tr.deg.$C\langle
x;\partial\rangle^{GL^{\partial}(m,F)}/C = n-m$.

Now to prove Theorem 1 it is enough to show that every
$x_1,...,x_n$ is $\delta$-algebraic over \\ $C\langle x,
\partial\rangle^{(GL^{\partial}(m,F),H)}.$ If one assumes
 that $\det [\delta^{\alpha^1}x,\delta^{\alpha^2}x,...,\delta^{\alpha^n}x]=0$ for all nonzero $\alpha^1,\alpha^2,...,\alpha^n$ from $W^n$, where $W$ stands for the set of whole numbers, then due to Proposition 2, applied to the differential field $(F\langle x, \partial\rangle,\delta )$, the system
$dx_1, dx_2,...,dx_n$ is linear dependent over $C$, because of
$\det [dx,d^2x,...,d^nx]=0$. So there is nontrivial system
$c_1,...,c_n$ elements of $C$ for which $\sum_{i=1}^n c_idx_i=
\sum_{i=1}^n c_idt\cdot \delta x_i= dt \cdot\sum_{i=1}^n c_i\delta
x_i=  dt \cdot\sum_{i=1}^n c_i\Phi^{\partial}\langle
x\rangle^{-1}\partial x_i =dt \cdot (\Phi^{\partial}\langle
x\rangle^{-1}\partial\sum_{i=1}^n c_i x_i )= 0$, where $dt$ stands
for row vector $(dt_1,dt_2,...,dt_m)$ and $\cdot$ for the dot
product. The last equality implies that $\sum_{i=1}^n c_i x_i =
0$, which can occur if and only if  $c_1=c_2=...=c_n=0$. This
contradiction shows that one can find nonzero
$\alpha^1,\alpha^2,...,\alpha^n$ from $W^n$ for which $\det
[\delta^{\alpha^1}x,\delta^{\alpha^2}x,...,\delta^{\alpha^n}x]\neq
0$. So now for any nonzero $\alpha \in W^n$ one can consider the
following differential equation in $y$:
$$\det [\delta^{\alpha^1}x,\delta^{\alpha^2}x,...,\delta^{\alpha^n}x]\det
[\delta^{\alpha^1}\overline{x},\delta^{\alpha^2}\overline{x},...,\delta^{\alpha^n}\overline{x}
,\delta^{\alpha}\overline{x}]=0$$,
where $\overline{x}=(x_1,x_2,...,x_n,y)$, $\delta^{\alpha}=\delta^{(\alpha_1,\alpha_2,...,\alpha_m)}$
stands for
$\delta_1^{\alpha_1}\delta_2^{\alpha_2}...\delta_m^{\alpha_m}.$
All coefficients of this differential equation belong to $C\langle
x, \partial\rangle^{(GL^{\partial}(m,F),H)}$ and $y= x_i$ is a
solution for this linear differential equation at any
$i=\overline{1,n}$. It implies that indeed
$\delta$-tr.deg.$C\langle x,
\partial\rangle^{(GL^{\partial}(m,F),H)}/C= n-m$.

The following result says that one can obtain a system of
generators of $(C\langle x\rangle^{(GL^{\partial}(m,F),H)}, \delta
)$ from the given system of generators of $(C\langle
x,\partial\rangle^{H},\partial )$.

{\bf Theorem 2.} \textit{If $(C\langle
x,\partial\rangle^{H},\partial )$ as a $\partial$-differential
field over $C$ is generated by a system
$(\varphi^{\partial}_{i}\langle x\rangle)_{i=\overline{1,l}}$ then
$\delta$-differential field $(C\langle x,
\partial\rangle^{(GL^{\partial}(m,F),H)}, \delta )$ is generated over $C$
by the system $(\varphi^{\delta}_{i}\langle
x\rangle)_{i=\overline{1,l}}$.}

{\bf Proof.} Let an irreducible $\frac{P^{\partial}\{
x\}}{Q^{\partial}\{ x\}}\in C\langle x,\partial\rangle^{H}$ be
$GL^{\partial}(m,F)$ -invariant. It means that for any $g\in
GL^{\partial}(m,F)$ one has the equality
$$P^{g^{-1}\partial}\{ x\}Q^{\partial}\{x\} =P^{\partial}\{ x\}Q^{g^{-1}\partial}\{ x\}$$
Therefore $P^{g^{-1}\partial}\{ x\} = P^{\partial}\{x\} \chi^{\partial}\langle g \rangle$.
 The function $\chi^{\partial}\langle T \rangle$( a "character" of $GL^{\partial}(m,F)$)
  has the following property
$$\chi^{\partial}\langle g_1g_2 \rangle =
\chi^{\partial}\langle g_1 \rangle\chi^{g^{-1}_1\partial}\langle
g_2 \rangle$$, for any $g_1 \in GL^{\partial}(m,F)$ and $g_2 \in
GL^{g^{-1}_1\partial}(m,F)$. But due to (2) one has $g_2=
g^{-1}_1g$ for some $g\in GL^{\partial}(m,F)$ therefore
$$\chi^{\partial}\langle g\rangle =\chi^{\partial}\langle g_1
\rangle\chi^{g^{-1}_1\partial}\langle g^{-1}_1g \rangle$$, for any $g_1,g \in GL^{\partial}(m,F)$.
It implies that
$$\chi^{\partial}\langle T\rangle =\chi^{\partial}\langle S \rangle\chi^{S^{-1}\partial}\langle S^{-1}T \rangle$$, for any $ T=(t_{ij})_{i,j=
\overline{1,m}}, S=(s_{ij})_{i,j= \overline{1,m}}$,  for which $
\partial_{k}t_{ij}=
 \partial_{i}t_{kj}, \partial_{k}s_{ij}=
 \partial_{i}s_{kj}$  at  $i,j,k= \overline{1,m} $.
The last equality guarantees that the function
$\chi^{\partial}\langle T\rangle $ can not vanish. Therefore
$\frac{P^{\partial}\{x\}}{Q^{\partial}\{ x\}}=\frac{P^{\delta}\{
x\}}{Q^{\delta}\{ x\}}$. This is the end of proof Theorem 2.

Let us assume that $(C\langle x,\partial \rangle^H, \partial)=
C\langle\varphi^{\partial}_1\langle x\rangle,
\varphi^{\partial}_2\langle
x\rangle,...,\varphi^{\partial}_{l}\langle x\rangle,\partial
\rangle.$ As far as all components of the matrix
$\Phi^{\partial}\langle x\rangle $ belong to   $C\langle
x,\partial \rangle^H$ it can be represented in the form
$$\Phi^{\partial}\langle x\rangle=\overline{\Phi}^{\partial}\langle \varphi^{\partial}_1\langle x\rangle, \varphi^{\partial}_2\langle x\rangle,...,\varphi^{\partial}_{l}\langle x\rangle\rangle =(\overline{\phi}^{\partial}_{ij}\langle \varphi^{\partial}_1\langle x\rangle, \varphi^{\partial}_2\langle x\rangle,...,\varphi^{\partial}_{l}\langle x\rangle\rangle )_{i,j=\overline{1,m}}$$, where $\overline{\phi}^{\partial}_{ij}\langle t_1,t_2,....,t_{l}\rangle \in C\langle t_1,t_2,....,t_{l} ,\partial \rangle$.
Therefore due to $\Phi^{\delta}\langle x\rangle= E_m$ one has
$$\overline{\Phi}^{\delta}\langle \varphi^{\delta}_1\langle
x\rangle, \varphi^{\delta}_2\langle
x\rangle,...,\varphi^{\delta}_{l}\langle x\rangle\rangle =E_m. $$

{\bf Remark 4.} The equality $$\chi^{\partial}\langle g\rangle
=\chi^{\partial}\langle g_1 \rangle\chi^{g^{-1}_1\partial}\langle
g^{-1}_1g \rangle$$, for any $g_1,g \in GL^{\partial}(m,F)$
resembles the property of character of the group $GL(m,F)$.
Therefore $\chi^{\partial}\langle T\rangle$ for which the above
equality is valid can be considered as a character of the groupoid
$GL^{\partial}(m,F)$. Description all such characters is an
interesting problem. Of course, $\chi^{\partial}\langle g\rangle=
\det(g)^k$, where $k$ is any integer number, are examples of such
characters. Are they all possible differential rational characters
of $GL^{\partial}(m,F)$?

{\bf Theorem 3.} \textit{Any $\delta$-differential polynomial
relation over $C$ of the system
 $\varphi^{\delta}_1\langle x\rangle, \varphi^{\delta}_2\langle x\rangle,...,\varphi^{\delta}_{l}\langle x\rangle $ is a
consequence of $\partial$-differential polynomial relations of the
system $\varphi^{\partial}_1\langle x\rangle,
\varphi^{\partial}_2\langle
x\rangle,...,\varphi^{\partial}_{l}\langle x\rangle $ over $C$ and
the relations $\overline{\Phi}^{\delta}\langle
\varphi^{\delta}_1\langle x\rangle, \varphi^{\delta}_2\langle
x\rangle,...,\varphi^{\delta}_{l}\langle x\rangle\rangle =E_m.$}

{\bf Proof.} Let $\psi^{\delta}\{\varphi^{\delta}_1\langle
x\rangle, \varphi^{\delta}_2\langle
x\rangle,...,\varphi^{\delta}_{l}\langle x\rangle \} = 0$, where
$\psi^{\partial}\{ t_1,t_2,...,t_{l} \} \in C\{
t_1,t_2,...,t_{l},\partial\}$.

If $\psi^{\partial}\{\varphi^{\partial}_1\langle x\rangle,
\varphi^{\partial}_2\langle
x\rangle,...,\varphi^{\partial}_{l}\langle x\rangle \} = 0$ then
it means that the above relation ($\psi^{\delta}\{
t_1,t_2,...,t_{l} \}$) of the system $\varphi^{\delta}_1\langle
x\rangle, \varphi^{\delta}_2\langle
x\rangle,...,\varphi^{\delta}_{l}\langle x\rangle $ is a
consequence of the relation ($\psi^{\partial}\{ t_1,t_2,...,t_{l}
\}$) of the system $\varphi^{\partial}_1\langle x\rangle,
\varphi^{\partial}_2\langle
x\rangle,...,\varphi^{\partial}_{l}\langle x\rangle $ i.e. it is
obtained by substitution $\delta$ for $\partial$ in
$\psi^{\partial}\{ t_1,t_2,...,t_{l} \}$.

If $\psi^{\partial}\{\varphi^{\partial}_1\langle x\rangle,
\varphi^{\partial}_2\langle
x\rangle,...,\varphi^{\partial}_{l}\langle x\rangle \} \neq 0$
then consider
$\psi^{T^{-1}\partial}\{\varphi^{T^{-1}\partial}_1\langle
x\rangle, \varphi^{T^{-1}\partial}_2\langle
x\rangle,...,\varphi^{T^{-1}\partial}_{l}\langle x\rangle \}$ as a
$\partial$-differential rational function in  variables
$T=(t_{ij})_{i,j=\overline{1,m}}$, where $\partial_kt_{ij}
=\partial_it_{kj}$ for any $i,j,k=\overline{1,m}$ over $C\langle
x,\partial\rangle$. Let $\frac{a^{\partial}_x\{
T\}}{b^{\partial}_x\{ T\}}$ be its irreducible representation and
the leading coefficient (with respect to some linear order) of
$b^{\partial}_x\{T\}$ be one. We show that in this case all
coefficients of $a^{\partial}_x\{ T\}, b^{\partial}_x\{T\}$ belong
to $C\langle x,\partial\rangle^H$.

Indeed, first of all
$\psi^{T^{-1}\partial}\{\varphi^{T^{-1}\partial}_1\langle
x\rangle, \varphi^{T^{-1}\partial}_2\langle
x\rangle,...,\varphi^{T^{-1}\partial}_{l}\langle x\rangle \}$, as
a differential rational function in $x$, is $H$- invariant
function, as much as $\varphi^{\partial}_i\langle x\rangle \in
C\langle x,\partial\rangle^H$. This $H$-invariantness implies that
$$a^{\partial}_x\{T\}b^{\partial}_{hx+h_0}\{ T\}= b^{\partial}_x\{T\}a^{\partial}_{hx+h_0}\{T\}$$ for any $(h,h_0)\in H$. Therefore
$b^{\partial}_{hx+h_0}\{ T\}= \chi^{\partial}\langle
x,(h,h_0)\rangle b^{\partial}_x\{T\}$. But comparision of the
leading terms of both sides implies that in reality
$\chi^{\partial}\langle x,(h,h_0)\rangle= 1$ which in its turn
implies that all coefficients of $b^{\partial}_x\{T\}$ (as well as
$a^{\partial}_x\{T\}$) are $H$- invariant.

Therefore all coefficients of $a^{\partial}_x\{T\}$,
$b^{\partial}_x\{T\}$ can be considered as
${\partial}$-differential rational functions in
$\varphi^{\partial}_1\langle x\rangle, \varphi^{\partial}_2\langle
x\rangle,...,\varphi^{\partial}_{l}\langle x\rangle$, for example,
let $b^{\partial}_x\{T\}=
\overline{b}_{\varphi^{\partial}_1\langle x\rangle,
\varphi^{\partial}_2\langle
x\rangle,...,\varphi^{\partial}_{l}\langle x\rangle}\{T\}$. Now
represent the numerator $a^{\partial}_x\{T\} $ as a
${\partial}$-differential polynomial function in $t_{ij}-
\overline{\phi}_{ij}^{\partial}\langle \varphi^{\partial}_1\langle
x\rangle, \varphi^{\partial}_2\langle
x\rangle,...,\varphi^{\partial}_{l}\langle x\rangle\rangle $,where
$i,j=\overline{1,m}$, for example, let
$a^{\partial}_x\{T\}=\overline{a}^{\partial}_{\varphi^{\partial}_1\langle
x\rangle, \varphi^{\partial}_2\langle
x\rangle,...,\varphi^{\partial}_{l}\langle x\rangle}\{T-
\overline{\phi}^{\partial}\langle \varphi^{\partial}_1\langle
x\rangle, \varphi^{\partial}_2\langle
x\rangle,...,\varphi^{\partial}_{l}\langle x\rangle \rangle\}$ .
As such polynomial its constant term is zero because of
$\psi^{\delta}\{\varphi^{\delta}_1\langle x\rangle,
\varphi^{\delta}_2\langle x\rangle,...,\varphi^{\delta}_{l}\langle
x\rangle \} = 0$. So
$$\psi^{T^{-1}{\partial}}\{\varphi^{T^{-1}{\partial}}_1\langle
x\rangle, \varphi^{T^{-1}{\partial}}_2\langle
x\rangle,...,\varphi^{T^{-1}{\partial}}_{l}\langle x\rangle \}
=\frac{\overline{ a}^{\partial}_{\varphi^{\partial}_1\langle
x\rangle, \varphi^{\partial}_2\langle
x\rangle,...,\varphi^{\partial}_{l}\langle
x\rangle}\{T-\overline{\phi}^{\partial}\langle
\varphi^{\partial}_1\langle x\rangle, \varphi^{\partial}_2\langle
x\rangle,...,\varphi^{\partial}_{l}\langle x\rangle \rangle
\}}{\overline{b}^{\partial}_{\varphi^{\partial}_1\langle x\rangle,
\varphi^{\partial}_2\langle
x\rangle,...,\varphi^{\partial}_{l}\langle x\rangle}\{T\}}.$$

Substitution $T= E_m$ implies that
$$\psi^{{\partial}}\{\varphi^{{\partial}}_1\langle x\rangle, \varphi^{{\partial}}_2\langle x\rangle,...,\varphi^{{\partial}}_{l}\langle x\rangle \} =\frac{\overline{
a}^{\partial}_{\varphi^{\partial}_1\langle x\rangle,
\varphi^{\partial}_2\langle
x\rangle,...,\varphi^{\partial}_{l}\langle
x\rangle}\{E_m-\overline{\phi}^{\partial}\langle
\varphi^{\partial}_1\langle x\rangle, \varphi^{\partial}_2\langle
x\rangle,...,\varphi^{\partial}_{l}\langle x\rangle \rangle
\}}{\overline{b}^{\partial}_{\varphi^{\partial}_1\langle x\rangle,
\varphi^{\partial}_2\langle
x\rangle,...,\varphi^{\partial}_{l}\langle x\rangle}\{E_m\}}.$$

Now consider the following $\delta$-differential rational function
over $C$:
$$\overline{\psi}^{\delta}\langle t_1,t_2,...,t_m \rangle=
\psi^{\delta}\{ t_1,t_2,...,t_{l}\}-\frac{\overline{
a}^{\delta}_{t_1,
t_2,...,t_{l}}\{E_m-\overline{\phi}^{\delta}\langle
t_1,t_2,...,t_{l} \rangle
\}}{\overline{b}^{\delta}_{t_1,t_2,...,t_{l}}\{E_m\}}.$$ For this
function one has
$$\overline{\psi}^{\delta}\langle \varphi^{\delta}_1\langle x\rangle, \varphi^{\delta}_2\langle x\rangle,...,\varphi^{\delta}_{l}\langle x\rangle \rangle= 0 \quad
\mbox{as well as} \quad \overline{\psi}^{\partial}\langle
\varphi^{\partial}_1\langle x\rangle, \varphi^{\partial}_2\langle
x\rangle,...,\varphi^{\partial}_{l}\langle x\rangle\rangle = 0.$$
But once again the last equality(relation) means that it is a
consequence of relations of the system
$\varphi^{\partial}_1\langle x\rangle, \varphi^{\partial}_2\langle
x\rangle,...,\varphi^{\partial}_{l}\langle x\rangle $. This is the
end of proof of Theorem 3.

{\bf  4. On $H$- invariants.}

The following result provides a method to find generators of the
differential field $C\langle x,\partial\rangle^H$ over $C$. Let
$\alpha^1, \alpha^2,...,\alpha^n$ be any different nonzero
elements of $W^n$.

For different classical subgroups $H$ of Affine group the field
$C\langle x,\partial\rangle^{H}$ is investigated in [7] in the
case of $m=1$. Our main concern here will be the case of $m\geq 1$
and arbitrary subgroup $H$ of affine group.

 {\bf Theorem 4.}
\textit{The equality $C\langle x,\partial\rangle^{GL(n,C)\propto
C^n}(x,\partial^{\alpha^1}x,\partial^{\alpha^2}x,...,\partial^{\alpha^n}x)=C\langle
x,\partial\rangle$ is valid and moreover the system consisting of
components of
$x,\partial^{\alpha^1}x,\partial^{\alpha^2}x,...,\partial^{\alpha^n}x$
is algebraic independent over  $C\langle
x,\partial\rangle^{GL(n,C)\propto C^n}$.}

{\bf Proof.} For any nonzero $\alpha \in W^n$ consider the
differential equation in $y$:
$$\det [\partial^{\alpha^1}x,\partial^{\alpha^2}x,...,\partial^{\alpha^n}x]\det
[\partial^{\alpha^1}\overline{x},\partial^{\alpha^2}\overline{x},...,\partial^{\alpha^n}
\overline{x},\partial^{\alpha}\overline{x}]=0$$, where
 $\overline{x}$ stands for $(x_1, x_2,..., x_n,y)$.
 All
 coefficients of this differential equation are in $$C\langle
x,\partial\rangle^{GL(n,C)\propto C^n}$$ and $y=x_i$ is a solution
whenever $i=1,2,...,n$. Therefore $$C\langle
x,\partial\rangle^{GL(n,C)\propto
C^n}(x,\partial^{\alpha^1}x,\partial^{\alpha^2}x,...,\partial^{\alpha^n}x)=C\langle
x,\partial\rangle$$

To prove algebraic independence of the system
$x,\partial^{\alpha^1}x,\partial^{\alpha^2}x,...,\partial^{\alpha^n}x$
over  $C\langle x,\partial\rangle^{GL(n,C)\propto C^n}$ it is
enough to show algebraic independence of the system
$\partial^{\alpha^1}x,\partial^{\alpha^2}x,...,\partial^{\alpha^n}x$
over  $C\langle x,\partial\rangle^{GL(n,C)\propto C^n}$.

Let $P[z^1,z^2,...,z^n]$, where $z^i= (z^i_1,z^i_2,...,z^i_n)$, be
a nonzero polynomial over  $C\langle
x,\partial\rangle^{GL(n,C)\propto C^n}$ such that
$P[\partial^{\alpha^1}x,\partial^{\alpha^2}x,,...,\partial^{\alpha^n}x]=
0$.  Assume, for example, at least one of \\$z^{i}_n$, where
$i=\overline{1,n}$, occurs in $P[z^1,z^2,...,z^n]$ and
$$P[z^1,z^2,...,z^n]= \sum_{\beta}(z^1_n)^{\beta_1}(z^2_n)^{\beta_2}...(z^n_n)^{\beta_n}P_{\beta}[\overline{z^1},\overline{z^2},...,\overline{z^n}]$$, where
$P_{\beta}[\overline{z^1},\overline{z^2},...,\overline{z^n}]$ are
polynomials over $C\langle x,\partial\rangle^{GL(n,C)\propto C^n}$
in  $\overline{z^i}=(z^i_1,z^i_2,...,z^i_{n-1})$,
$i=\overline{1,n}$.

Consider $h \in GL(n,C)$ which's $i$-th column is of the form
$(0,...,0,1,0,...,0,c_i)$, where $i=\overline{1,n-1}$ and its
$n$-th column is $(0,...,0,c_n)$. For such $h$ one has
$\overline{hx}= \overline{x}$.
 So far as the coefficients of $P[z^1,z^2,...,z^n]$ are $GL(n,C)\propto C^n$- invariant, substitution $hx$ for $x$ into \\
$P[\partial^{\alpha^1}x,\partial^{\alpha^2}x,,...,\partial^{\alpha^n}x]=
0$
 implies that
$$\sum_{\beta}(\sum_{i=1}^nc_i\partial^{\alpha^1}x_i)^{\beta_1}(\sum_{i=1}^nc_i\partial^{\alpha^2}x_i)^{\beta_2}...(\sum_{i=1}^nc_i\partial^{\alpha^n}x_i)^{\beta_n}P_{\beta}[\partial^{\alpha^1}\overline{x},\partial^{\alpha^2}\overline{x},...,\partial^{\alpha^n}\overline{x}]=0.$$  Therefore due to the second assumption on $(F,\partial)$ for variables $y_1,y_2,...,y_n$ one has
\begin{eqnarray}\begin{array}{c}
\sum_{\beta}(\sum_{i=1}^ny_i\partial^{\alpha^1}x_i)^{\beta_1}(\sum_{i=1}^ny_i\partial^{\alpha^2}x_i)^{\beta_2}...(\sum_{i=1}^ny_i\partial^{\alpha^n}x_i)^{\beta_n}P_{\beta}[\partial^{\alpha^1}\overline{x},\partial^{\alpha^2}\overline{x},...,\partial^{\alpha^n}\overline{x}]=0
\end{array}\end{eqnarray}

 Now consider the ring $C\langle x,\partial\rangle[y_1,y_2,...,y_n]$ with respect to the differential operators\\ $\overline{\partial}_1= \frac{\partial}{\partial y_1},\overline{\partial}_2= \frac{\partial}{\partial y_2},
...,\overline{\partial}_n= \frac{\partial}{\partial y_n}$. It is
clear that its constant ring is $C\langle x,\partial\rangle$ i.e.
$$C\langle x,d\rangle= \{a\in C\langle
x,\partial\rangle[y_1,y_2,...,y_n]:
\overline{\partial}_1a=\overline{\partial}_2a=...=\overline{\partial}_na=0
\}.$$ Introduce new differential operators
$\overline{\overline{\partial}}=
\sum_{j=1}^nf^{\partial}_{ij}\langle x
\rangle\overline{\partial}_j$, where
$i=\overline{1,n}$,$$(f^{\partial}_{ij}\langle x
\rangle)_{i,j=\overline{1,n}}=[\partial^{\alpha^1}x,\partial^{\alpha^2}x,...,\partial^{\alpha^n}x
]^{-1}$$ The following are evident:

 a) The constant ring of $(C\langle x,\partial\rangle [y_1,y_2,...,y_n],\overline{\overline{\partial}})$, where $\overline{\overline{\partial}}=(\overline{\overline{\partial}}_1,\overline{\overline{\partial}}_2,...,\overline{\overline{\partial}}_n)$ is the same $C\langle x,\partial\rangle$,

b)
$\overline{\overline{\partial}}_j(\sum_{i=1}^ny_i\partial^{\alpha^k}x_i)$
is equal to  $0$ whenever
 $j \neq k$ and it is equal to 1 if $j= k$, where $j,k=\overline{1,n}$ .

Now if one assumes that $\beta^0= (\beta^0_1,...,\beta^0_n)$ is a
such one for which $$\vert \beta^0 \vert = \mbox{ max}\{
\vert\beta \vert:
P_{\beta}[\partial^{\alpha^1}\overline{x},\partial^{\alpha^2}\overline{x},...,\partial^{\alpha^n}\overline{x}]\neq
0 \}$$ and applies
$\overline{\overline{\partial}}_1^{\beta^0_1}\overline{\overline{\partial}}_2^{\beta^0_2}...\overline{\overline{\partial}}_n^{\beta^0_n}$
to equality (4) he comes to a contradiction
$P_{\beta^0}[\partial^{\alpha^1}\overline{x},\partial^{\alpha^2}\overline{x},...,\partial^{\alpha^n}\overline{x}]=
0$. This is the end of proof Theorem 4.

So due to Theorem 4 $$C\langle x, \partial\rangle^H= C\langle x,
\partial\rangle^{GL(n,C)\propto
C^n}(x,\partial^{\alpha^1}x,...,\partial^{\alpha^n}x)^H$$ and the
system $x,\partial^{\alpha^1}x,...,\partial^{\alpha^n}x$ is
algebraic independent over $C\langle x,
\partial\rangle^{GL(n,C)\propto C^n}$. Note that every element of
the field $C\langle x, \partial\rangle^{GL(n,C)\propto C^n}$ is a
fixed element for the group $H$. Therefore if one wants to have a
system of differential generators of $(C\langle x,
\partial\rangle^H, \partial)$ over $C$ he can do the following:

1. Find any system of generators (over C) of the differential
field
$$(C\langle x, \partial\rangle^{GL(n,C)\propto C^n},\partial )$$

2. Find any system of ordinary algebraic generators of the field
$$C\langle x, \partial\rangle^{GL(n,C)\propto
C^n}(z^1,z^2,...,z^{n+1})^H$$, where
$z^i=(z^i_1,z^i_2,...,z^i_{n})$ , $i=\overline{1,n+1}$, and the
action of $H$ is defined as:
$$((h,h_0),(z^1,z^2,...,z^{n+1}))\rightarrow
(hz^1+h_0,hz^2,...,hz^{n+1})$$  For example, let it be
$\varphi_1(z^1,z^2,...,z^{n+1}),
\varphi_2(z^1,z^2,...,z^{n+1}),...,\varphi_k(z^1,z^2,...,z^{n+1})$.

Then the union of the system of generators of $(C\langle x,
\partial\rangle^{GL(n,C)\propto C^n},\partial )$ with $$\{
\varphi_1(x,\partial^{\alpha^1}x,...,\partial^{\alpha^n}x),
\varphi_2(x,\partial^{\alpha^1}x,...,\partial^{\alpha^n}x),...,\varphi_k(x,\partial^{\alpha^1}x,...,\partial^{\alpha^n}x)\}$$
can be taken as a system of generators of the differential field
$(C\langle x, \partial\rangle^H, \partial)$ over $C$.

In the case of $m=1$ to find a system of generators of the field
$C\langle x, \partial\rangle^{GL(n,C)\propto
C^n}(z^1,z^2,...,z^{n+1})^H$ it was enough to find generators of
$C(z^1,z^2,...,z^{n+1})^H$ due to the fact that the differential
field $(C\langle x,
\partial\rangle^{GL(n,C)\propto C^n},\partial )$ has a
$\partial$-algebraic independent system of generators over $C$. It
seems that if $m>1$ this fact is not true for $(C\langle x,
\partial\rangle^{GL(n,C)\propto C^n},\partial )$ anymore.

{\bf Remark 5.} At the end I would like to note that a different
approach can be done to the equivalence problem of surfaces by the
use of rational differential forms. Definition of the ordinary
high order differentials of (for example, real) functions of $m$
variables can be changed slightly in such a way that not only the
first differential but also all high order differentials will have
invariant form with respect to change of variables [8]. One can
use
 it to introduce the field of differential rational forms $R\langle x,
d\rangle$, where $x=(x_1,...,x_n)$ is assumed to be variable
$m$-parametric surface. Moreover this field is an ordinary
differential field with respect to $d$. Due to the invariant
property of high order differentials with respect to change of
variables one have to consider only $R\langle x, d\rangle^H$,
where $H$ is a given motion group of $R^n$. The geometric meaning
(or interpretation) of such differential rational forms are not
clear but nevertheless one can use results from [2] to find a
system of generators of the differential field  ($R\langle
x,d\rangle^H, d).$ \newpage

\begin{center}{References}\end{center}

1. A.Weinstein, Groupoids: Unifying Internal and External
Symmetry. Notices of the AMS, Volume 43, Number 7, 744-752.

2. Ural Bekbaev, On the field of differential rational invariants
of a subgroup of Affine group(Ordinary differential case).
arXiv:math: AG/0608479

3. Bekbaev U., Differential rational invariants of Hypersurfaces
relative to Affine Group. In  Prociding Pengintegrasion Teknologi
dalam Sains Matematik, ed. How Guang Aun, Leong Fook, Ong Boon
Hua, Quah Soon Hoe, Safian Uda, Zarita Zainiddin, pp. 58-64, Pusat
Pengajian Sains Matemarik, USM, 1999.

4. Bekbaev U., On differential rational invariants of finite
subgroups of Affine group. Bulletin of Malaysian Mathematical
Society, 2005, Volume 28, N1, pp. 55-60.

5. Bekbaev U.D., An algebraic approach to invariants of surfaces.
Proceedings of IPTA Research \& Development Exposition 2003. 9-12
October 2003. Putra World Trade Center, Kuala Lumpur. Vol. 5:
Science and Engineering, ed. Y.A. Khalid  et al.,pp. 299-306.
University Putra Malaysia Press, Serdang, Selangor, Malaysia.

6. E.R.Kolchin, " Differential Algebra and Algebraic Groups",
Academic Press, New York, 1973.

7. Dj.Khadjiev, "Application of Invariant  Theory to Differential
Geometry of curves", FAN, Tashkent, 1988(Russian).

8. U.D. Bekbaev, High order invariant differentials of functions
in local coordinates. Dep.v VINITI, 25.07.90, N4225-B90 (Russian).

\end{document}